\documentclass[12pt]{article}
\usepackage{amssymb,amsthm,amscd,array}

\let\ssection=\section
\renewcommand{\section}{\setcounter{equation}{0}\ssection}

\newcommand{\bbR}{\mathbb{R}}

\newcommand{\aff}{\mathrm{aff}}

\newcommand{\tC}{\widetilde{C}}

\newcommand{\cD}{{\mathcal{D}}}

\newcommand{\cE}{{\mathcal{E}}}

\newcommand{\rD}{\mathrm{D}}
\newcommand{\rE}{\mathrm{E}}

\newcommand{\cF}{{\mathcal{F}}}

\newcommand{\gr}{{\mathrm{gr}}}
\newcommand{\rg}{\mathrm{g}}

\newcommand{\cI}{{\mathcal{I}}}

\newcommand{\cQ}{{\mathcal{Q}}}

\newcommand{\cS}{{\mathcal{S}}}

\newcommand{\SL}{\mathrm{SL}}
\newcommand{\Sl}{\mathrm{sl}}
\newcommand{\SO}{\mathrm{SO}}

\newcommand{\half}{\frac{1}{2}}

\begin{document}

\baselineskip=18pt

%\textwidth=16truecm
%\textheight=24truecm
%\hoffset=-1.5truecm
%\voffset=-2.5truecm

\def\a{\alpha}
\def\b{\beta}
\def\d{\delta}
\def\g{\gamma}
\def\om{\omega}
\def\r{\rho}
\def\s{\sigma}
\def\vfi{\varphi}
\def\l{\lambda}
\def\m{\mu}
\def\implies{\Rightarrow}

\oddsidemargin .1truein
\newtheorem{thm}{Theorem}[section]
\newtheorem{lem}[thm]{Lemma}
\newtheorem{cor}[thm]{Corollary}
\newtheorem{pro}[thm]{Proposition}
\newtheorem{ex}[thm]{Example}
\newtheorem{rmk}[thm]{Remark}
\newtheorem{defi}[thm]{Definition}
%\newremark{ex}[thm]{Example}

%\newtheorem{thm}{Theorem}
%\newtheorem{lem}{Lemma}
%\newtheorem{cor}{Corollary}
%\newtheorem{prop}[thm]{Proposition}
%\newtheorem{definition}{Definition}

\title{Projectively equivariant quantization \\and symbol calculus:\\
noncommutative hypergeometric functions}

\author{C.~Duval\footnote{mailto:duval@cpt.univ-mrs.fr}\\
{\small Universit\'e de la M\'editerran\'ee and CPT-CNRS}
\and
V.~Ovsienko\footnote{mailto:ovsienko@cpt.univ-mrs.fr}\\
{\small CNRS, Centre de Physique Th\'eorique}
\thanks{CPT-CNRS, Luminy Case 907,
F--13288 Marseille, Cedex 9, FRANCE.
}
}

\date{}

\maketitle

\thispagestyle{empty}

\begin{abstract}
We extend projectively equivariant quantization and symbol calculus to symbols of
pseudo-differential operators. An explicit expression in terms of hypergeometric
functions with noncommutative arguments is given. Some examples are worked out, one of
them yielding a quantum length element on~$S^3$.
\end{abstract}

\vskip1cm
\noindent
\textbf{Keywords:} Quantization, projective structures, hypergeometric functions.

\newpage

%%%%%%%%%%%%%%%%%%%%%%%%%%%%%%%%%%%%%%%%%%%%%%%%%%%%%%%%%%%%%%%%%%%%%%%%%%%%%%%%%%%%
%%%%%%%%%%%%%%%%%%%%%%%%%%%%%%%%%%%%%%%%%%%%%%%%%%%%%%%%%%%%%%%%%%%%%%%%%%%%%%%%%%%%
\section{Introduction}
%%%%%%%%%%%%%%%%%%%%%%%%%%%%%%%%%%%%%%%%%%%%%%%%%%%%%%%%%%%%%%%%%%%%%%%%%%%%%%%%%%%%
%%%%%%%%%%%%%%%%%%%%%%%%%%%%%%%%%%%%%%%%%%%%%%%%%%%%%%%%%%%%%%%%%%%%%%%%%%%%%%%%%%%%

Let $M$ be a smooth manifold and $\cS(M)$ the space of smooth functions on
$T^*M$, polynomial on the fibers; the latter is usually called the space of symbols of
differential operators. Let us furthermore assume that $M$ is endowed with an action of a
Lie group $G$. The aim of equivariant quantization  \cite{LO,DO2,DLO} (see also
\cite{CMZ}, and \cite{B1,B2}) is to associate to each symbol a differential operator on
$M$ in such a way that this quantization map intertwines the $G$-action.

The existence and uniqueness of equivariant quantization  in the
case where~$M$ has a flat projective (resp. conformal) structure, i.e., when
$G=\SL(n+1,\bbR)$ with $n=\dim(M)$ (resp. $G=\SO(p+1,q+1)$ with $p+q=\dim(M)$) has
recently been proved in the above references.

More precisely, let $\cF_\l(M)$ stand for the space of (complex-valued) tensor densities
of degree $\l$ on $M$ and $\cD_{\l,\m}(M)$ for the space of linear differential operators
from $\cF_\l(M)$ to $\cF_\m(M)$. These spaces are naturally modules over the group of all
diffeomorphisms of $M$. The space of symbols corresponding to
$\cD_{\l,\m}(M)$ is therefore $\cS_\d(M)=\cS(M)\otimes\cF_\d(M)$ where $\d=\mu-\l$. There
is a filtration
$$
\cD^0_{\l,\m}\subset\cD^1_{\l,\m}\subset\cdots\subset\cD^k_{\l,\m}\subset\cdots
$$ 
and the associated module $\cS_\d(M)=\gr\left(\cD_{\l,\m}\right)$ is graded by the degree
of polynomials: 
$$
%%%\cS_\d=\bigoplus_{k=0}^\infty{\cS_{k,\d}}
\cS_\d=\cS_{0,\d}\oplus\cS_{1,\d}\oplus\cdots\oplus\cS_{k,\d}\oplus\cdots
$$ 
The problem of equivariant quantization is the quest for a
\textit{quantization map}:
\begin{equation}
\cQ_{\l,\m}:\cS_\d(M)\to\cD_{\l,\m}(M)
\label{Q}
\end{equation}
that commutes with the $G$-action.
In other words, it amounts to an identification of these two spaces which is canonical
with respect to the geometric structure on $M$.

The inverse of the quantization map:
\begin{equation}
\s_{\l,\m}=\left(\cQ_{\l,\m}\right)^{-1}
\label{sigma}
\end{equation}
is called the symbol map.

\goodbreak

In this Letter, we will restrict considerations to the projectively equivariant case.
Without loss of generality, we will assume $M=S^n$ endowed with its standard
$\SL(n+1,\bbR)$-action. The explicit formul{\ae} for the maps (\ref{Q}) and (\ref{sigma})
can be found in
\cite{CMZ} for $n=1$ and in
\cite{LO} for
$\l=\m$ in any dimension.
Our purpose is to rewrite the expressions for $\cQ_{\l,\m}$ and $\s_{\l,\m}$ in a more
general way which, in particular, extends the quantization to a bigger class of symbols
of pseudo-differential operators.

%%%%%%%%%%%%%%%%%%%%%%%%%%%%%%%%%%%%%%%%%%%%%%%%%%%%%%%%%%%%%%%%%%%%%%%%%%%%%%%%%%%%
%%%%%%%%%%%%%%%%%%%%%%%%%%%%%%%%%%%%%%%%%%%%%%%%%%%%%%%%%%%%%%%%%%%%%%%%%%%%%%%%%%%%
\section{Projectively equivariant quantization map}
%%%%%%%%%%%%%%%%%%%%%%%%%%%%%%%%%%%%%%%%%%%%%%%%%%%%%%%%%%%%%%%%%%%%%%%%%%%%%%%%%%%%
%%%%%%%%%%%%%%%%%%%%%%%%%%%%%%%%%%%%%%%%%%%%%%%%%%%%%%%%%%%%%%%%%%%%%%%%%%%%%%%%%%%%

In terms of affine coordinates on $S^n$, the vector fields spanning the canonical
action of the Lie algebra $\Sl(n+1,\bbR)$ are as follows
$$
\frac{\partial}{\partial{}x^i},
\qquad
x^i\frac{\partial}{\partial{}x^j},
\qquad
x^ix^j\frac{\partial}{\partial{}x^j},
$$
with $i,j=1,\ldots,n$ (the Einstein summation convention is understood). 

Will will denote
by $\aff(n,\bbR)$ the affine subalgebra spanned by the first-order vector fields.
We will find it convenient to identify locally, in each affine chart, the spaces
$\cS_\d$ and
$\cD_{\l,\m}$ via the ``normal ordering'' isomorphism
\begin{equation}
\cI:P(x)^{i_1\ldots{}i_k}\,\xi_{i_1}\cdots\xi_{i_k}
\longmapsto
(-i\hbar)^k\,P(x)^{i_1\ldots{}i_k}\frac{\partial}{\partial{}x^{i_1}}
\cdots
\frac{\partial}{\partial{}x^{i_k}}
\label{I}
\end{equation}
which is already equivariant with respect to $\aff(n,\bbR)$. An equivalent means of
identification is provided by the Fourier transform
\begin{equation}
\left(\cI(P)\phi\right)(x)
=
\frac{1}{(2\pi\hbar)^{n/2}}
\int_{\bbR^{2n}}{\!\!\!e^{(i/\hbar)\langle\xi,x-y\rangle}P(y,\xi)\phi(y)\,dyd\xi}
\label{FT}
\end{equation}
where $P(y,\xi)=\sum_{k=0}^{\infty}{P(y)^{i_1\ldots{}i_k}\,\xi_{i_1}\cdots\xi_{i_k}}$ and
where $\phi$ is a compactly supported function (representing a
$\l$-density in the coordinate patch). 
This mapping extends to the space of pseudo-differential symbols (defined in the chosen
affine coordinate system).

The purpose of projectively equivariant quantization is to modify the map~$\cI$
in~(\ref{I}) in order to obtain an identification of $\cS_\d$ and $\cD_{\l,\m}$ that does
not depend upon a chosen affine coordinate system, and is, therefore, globally defined on
$S^n$.

Recall \cite{Wey} that the (locally defined) operators on $\cS_\d$, namely
\begin{equation}
\cE=
\xi_i\frac{\partial}{\partial{}\xi_i},
\qquad
\rD=\frac{\partial}{\partial{}x^i}\frac{\partial}{\partial{}\xi_i},
\label{ED}
\end{equation}
(where the $\xi_i$ are the coordinates dual to the $x^i$) commute with the
$\aff(n,\bbR)$-action on $T^*S^n$. The Euler operator,
$\cE$, is the degree operator on $\cS_\d=\bigoplus_{k=0}^\infty{\cS_{k,\d}}$ while the
divergence operator
$\rD$ lowers this degree by one.

Let us now recall (in a slightly more general context) the results obtained
in~\cite{LO,CMZ}. The $\SL(n+1,\bbR)$-equivariant quantization map (\ref{Q}) is given on
every homogeneous component by
\begin{equation}
\cQ_{\l,\m}\big|_{\cS_{k,\d}}
=
\sum_{m=0}^k{
C_m^k\,(i\hbar\rD)^m|_{\cS_{k,\d}}
}
\label{Qaff}
\end{equation}
where the constant coefficients $C_m^k$ are determined by the following
relation
\begin{equation}
C_{m+1}^k
=
\frac{k-m-1+(n+1)\l}{(m+1)(2k-m-2+(n+1)(1-\d))}\,
C_{m}^k
\label{relations}
\end{equation}
and the normalization condition: $C_0^k=1$. 

As to the projectively equivariant symbol map (\ref{sigma}), it retains the form
\begin{equation}
\s_{\l,\m}\big|_{\cS_{k,\d}}
=
\sum_{m=0}^k{
\tC_m^k\left(\frac{\rD}{i\hbar}\right)^m\Big|_{\cS_{k,\d}}
}
\label{sigmaaff}
\end{equation}
where the coefficients $\tC_m^k$ are such that
\begin{equation}
\tC_{m+1}^{k+1}
=
-\frac{k+(n+1)\l}{(m+1)(2k-m+(n+1)(1-\d))}\,
\tC_{m}^k
\label{sigmarelations}
\end{equation}
and, again, $\tC_0^k=1$ for all $k$.

\begin{rmk}
{\rm
The expressions (\ref{Qaff}) and (\ref{sigmaaff}) make sense if 
%%% \begin{equation}
$$
\d\neq1+\frac{\ell}{n+1}
$$
%%% \label{resonances}
%%% \end{equation}
with
$\ell=0,1,2,\ldots$ For these values of $\d$, the quantization and symbol maps do not
exist for generic $\l$ and $\m$; see \cite{Lec1} for a detailed classification.
}
\end{rmk}

In contradistinction with the operators $\cE$ and~$\rD$ defined in (\ref{ED}), the
quantiza\-tion map $\cQ_{\l,\m}$ and the symbol map
$\s_{\l,\m}$ are globally defined on $T^*S^n$, i.e., they are independent of the choice
of an affine coordinate system.

%%%%%%%%%%%%%%%%%%%%%%%%%%%%%%%%%%%%%%%%%%%%%%%%%%%%%%%%%%%%%%%%%%%%%%%%%%%%%%%%%%%%
%%%%%%%%%%%%%%%%%%%%%%%%%%%%%%%%%%%%%%%%%%%%%%%%%%%%%%%%%%%%%%%%%%%%%%%%%%%%%%%%%%%%
\section{Noncommutative hypergeometric function}
%%%%%%%%%%%%%%%%%%%%%%%%%%%%%%%%%%%%%%%%%%%%%%%%%%%%%%%%%%%%%%%%%%%%%%%%%%%%%%%%%%%%
%%%%%%%%%%%%%%%%%%%%%%%%%%%%%%%%%%%%%%%%%%%%%%%%%%%%%%%%%%%%%%%%%%%%%%%%%%%%%%%%%%%%

Our main purpose is to obtain an expression for $\cQ_{\l,\m}$ and $\s_{\l,\m}$ valid for
a larger class of symbols, namely for symbols of \textit{pseudo-differential} operators.
We will rewrite the formul{\ae} (\ref{Qaff}), (\ref{relations}) and (\ref{sigmaaff}),
(\ref{sigmarelations}) in terms of the $\aff(n,\bbR)$-invariant operators $\cE$ and $\rD$
in a form independent of the degree, $k$, of polynomials. 

It turns out that our quantization map (\ref{Q}) involves a certain hyper\-geometric
function; let us now recall this classical notion.
A hypergeometric function with $p+q$ parameters is defined (see, e.g., \cite{GKP}) as the
power series in $z$ given by
\begin{equation}
F\left(\left.
\begin{array}{c}
a_1,\ldots,a_p\\b_1,\ldots,b_q
\end{array}\right|z\right)
=
\sum_{m=0}^\infty
{
\frac{(a_1)_m\cdots(a_p)_m}{(b_1)_m\cdots(b_q)_m}\,\frac{z^m}{m!}
}
\label{F}
\end{equation}
with $(a)_m=a(a+1)\cdots(a+m-1)$. This hyper\-geometric function is called confluent if
$p=q=1$.

\begin{thm}\label{Qthm}
The projectively equivariant quantization map is of the form
\begin{equation}
\cQ_{\l,\m}=F\left(\left.
\begin{array}{c}
A_1,A_2\\
B_1,B_2
\end{array}
\right|Z\right)
\label{QPE}
\end{equation}
where the parameters
\begin{equation}
\begin{array}{ll}
A_1=\cE+(n+1)\l,&A_2=2\cE+(n+1)(1-\d)-1,\\[8pt]
B_1=\cE+\half(n+1)(1-\d)-\half,&B_2=\cE+\half(n+1)(1-\d),
\end{array}
\label{AB}
\end{equation}
are operator-valued, as well as the variable
\begin{equation}
Z=\frac{i\hbar\rD}{4}.
\label{Z}
\end{equation}
\end{thm}
\begin{proof}
Recall that for a hypergeometric function (\ref{F}), one has
$$
F\left(\left.
\begin{array}{c}
a_1,\ldots,a_p\\b_1,\ldots,b_q
\end{array}\right|z\right)
=
\sum_{m=0}^\infty
{
c_m\,z^m
}
$$
with
$$
\frac{c_{m+1}}{c_m\hfill}
=
\frac{1}{m+1}
\left[
\frac{(a_1+m)\cdots(a_p+m)}{(b_1+m)\cdots(b_q+m)}
\right].
$$
Let us replace $k-m$ by the degree operator $\cE$ in the coefficients $C_m^k$; the
expression (\ref{Qaff}) can be therefore rewritten as 
$
\cQ_{\l,\m}
=
\sum_{m=0}^k{
C_m(\cE)\,(i\hbar\rD)^m
}
$.
From the recursion relation (\ref{relations}), one readily obtains
$$
\frac{C_{m+1}(\cE)}{C_m(\cE)\hfill}
=
\frac{1}{4(m+1)}\left[
\frac{\left(\cE+(n+1)\l+m\right)\left(2\cE+(n+1)(1-\d)-1+m\right)}
{\left(\cE+\half(n+1)(1-\d)-\half+m\right)\left(\cE+\half(n+1)(1-\d)+m\right)}
\right]
$$
completing the proof.
\end{proof}

\begin{cor}
The quantization map is given by the series
\begin{equation}
\cQ_{\l,\m}
=
\sum_{m=0}^k{
C_m(\cE)\,(i\hbar\rD)^m
}
\label{QPEbis}
\end{equation}
where
\begin{equation}
C_m(\cE)
=
\frac{1}{m!}
\frac{(\cE+(n+1)\l)_m}{\left(2\cE+(n+1)(1-\d)+m-1\right)_m}.
\label{Coeffm}
\end{equation}
\end{cor}

\goodbreak

\begin{rmk}
{\rm
Let us stress that the operator-valued parameters (\ref{AB}) and the variable
(\ref{Z}) entering the expression (\ref{QPE}) do not commute. We have therefore chosen an
ordering that assigns the divergence operator $\rD$ to the right.
}
\end{rmk}

In the particular and most interesting case of half-densities (cf. \cite{DO2,DLO}), the
expression~(\ref{QPE}) takes a simpler form.
\begin{cor}
If $\l=\m=\half$, the quantization map (\ref{QPE}) reduces to the confluent
hypergeometric~function
\begin{equation}
\cQ_{\half,\half}=F\left(\left.
\begin{array}{r}
2\rE\\\rE
\end{array}\right|\frac{i\hbar\rD}{4}\right)
\label{QPEhalf}
\end{equation}
with the notation: $\rE=\cE+\half\,n$.
\end{cor}

It is a remarkable fact that the expression for inverse symbol map (\ref{sigma})
is much simpler. It is given by a confluent hypergeometric function for any $\l$ and
$\m$.

\begin{thm}
The projectively equivariant symbol map (\ref{sigma}) is given by
\begin{equation}
\s_{\l,\m}=F\left(\left.
\begin{array}{l}
\cE+(n+1)\l\\
2\cE+(n+1)(1-\d)
\end{array}
\right|-\frac{D}{i\hbar}\right).
\label{SPE}
\end{equation}
\label{sigmathm}
\end{thm}
\noindent
The proof is analogous to that of Theorem \ref{Qthm}.

It would be interesting to obtain expressions of the projectively equivariant quantization
and symbol maps as integral operators similar to (\ref{FT}).

%%%%%%%%%%%%%%%%%%%%%%%%%%%%%%%%%%%%%%%%%%%%%%%%%%%%%%%%%%%%%%%%%%%%%%%%%%%%%%%%%%%%
%%%%%%%%%%%%%%%%%%%%%%%%%%%%%%%%%%%%%%%%%%%%%%%%%%%%%%%%%%%%%%%%%%%%%%%%%%%%%%%%%%%%
\section{Some examples}
%%%%%%%%%%%%%%%%%%%%%%%%%%%%%%%%%%%%%%%%%%%%%%%%%%%%%%%%%%%%%%%%%%%%%%%%%%%%%%%%%%%%
%%%%%%%%%%%%%%%%%%%%%%%%%%%%%%%%%%%%%%%%%%%%%%%%%%%%%%%%%%%%%%%%%%%%%%%%%%%%%%%%%%%%

We wish to present, here, a few applications of the  projectively equivariant quantization
to some special Hamiltonians on $T^*S^n$. 

The first example deals with the geodesic flow. Denote by $\rg$ the standard round metric
on the unit $n$-sphere and by $H=\rg^{ij}\xi_i\xi_j$ the corresponding quadratic
Hamiltonian. In an affine coordinate system, it takes the following form 
\begin{equation}
H=
\left(1+\Vert{x}\Vert^2\right)
\left(\delta^{ij}+x^ix^j\right)
\xi_i\xi_j
\label{H}
\end{equation}
where $\Vert{x}\Vert^2=\delta_{ij}x^ix^j$ with $i,j=1,\ldots,n$.
Moreover, we will consider a family of such Hamiltonians belonging to
$\cS_\d$, namely $H_\d=H\,\sqrt{g^\d}$ where $g=\det(\rg_{ij})$. 

In order to provide explicit formul{\ae}, we need to recall the expression of the
covariant derivative of $\l$-densities, namely
$\nabla_i=\partial_i-\l\Gamma^j_{ij}$.

\begin{pro} The projectively equivariant quantization map (\ref{Q}) associates to~$H_\d$
the following differential operator
\begin{equation}
\cQ_{\l,\m}(H_\d)=
-\hbar^2\left(
\Delta + C_{\l,\m}\,R
\right)
\label{QH}
\end{equation}
where $\Delta=\rg^{ij}\nabla_i\nabla_j$ is the Laplace operator;
the constant coefficient~is
\begin{equation}
C_{\l,\m}=\frac{(n+1)^2\l(\m-1)}{(n-1)((1-\d)(n+1)+1)}
\label{C}
\end{equation}
and $R=n(n-1)$ is the scalar curvature of $S^n$.
\end{pro}
\begin{proof}
The quantum operator (\ref{QH}) is obtained, using (\ref{QPE})--(\ref{Z}), by a direct
computation. However, the formula (\ref{QH}) turns out to be a particular case of (5.4)
in~\cite{Bo} since the Levi-Civita connection is projectively flat.
\end{proof}

Another example is provided by the $\a$-th power $H^\a$ of the Hamiltonian $H$,
where~$\a\in\bbR$. We will only consider the case $\l=\m$ in the sequel.
\begin{pro}
For
\begin{equation}
\a=\frac{1-n}{4}
\label{alpha}
\end{equation}
one has $\cQ_{\l,\l}(H^\a)=H^\a$.
\end{pro}
\begin{proof}
Straightforward computation leads to
$$
D(H^\a)=2\a(4\a+n-1)H^{\a-1}(1+\Vert{x}\Vert^2)\langle\xi,x\rangle
$$
and (\ref{Qaff}) therefore yields the result.
\end{proof}

\goodbreak

We have just shown that the Fourier transform
(\ref{FT}) of $H^{(1-n)/4}$ is well-defined on $S^n$ and actually corresponds to the
projectively equivariant quantization of this pseudo-differential symbol.

If we want to deal with operators acting on a Hilbert space, we have to restrict now
considerations to the case $\l=\m=\half$.

For the $3$-sphere only, the above quantum Hamiltonian on $T^*S^n\setminus{}S^n$ is as
follows
\begin{equation}
\cQ_{\half,\half}(H^{-\half})=\frac{1}{\hbar}\,\frac{1}{\sqrt{-\Delta}}
\label{ds}
\end{equation}
and can be understood as a quantized ``length element'' in the sense of \cite{Co}.

\bigskip
\noindent
\textbf{Acknowledgements:}
We thank P.~Lecomte and E.~Mourre for numerous enlightening discussions.

\goodbreak

%%%\newpage

%%%%%%%%%%%%%%%%%%%%%%%%%%%%%%%%%%%%%%%%%%%%%%%%%%%%%%%%%%%%%%%%%%%%%%%%%%%%%%%%%%%%
%%%%%%%%%%%%%%%%%%%%%%%%%%%%%%%%%%%%%%%%%%%%%%%%%%%%%%%%%%%%%%%%%%%%%%%%%%%%%%%%%%%%

\end{document}